\documentclass{tac}

\usepackage{amssymb,amsmath}
\usepackage{tikz-cd}
\usepackage{commath}
\usepackage{empheq}

\usepackage{xy}
\input diagxy

\usepackage[colorlinks=true]{hyperref}
\hypersetup{allcolors=[rgb]{0.1,0.1,0.4}}

\author{Daniel Graves}

\address{School of Mathematics, University of Leeds, Woodhouse, Leeds, LS2 9JT, UK
}

\title{Composing PROBs}

\copyrightyear{2022}


\eaddress{dan.graves92@gmail.com}

\renewcommand{\ul}{\underline}
\renewcommand{\phi}{\varphi}
\newcommand{\llb}{\left\lbrace}
\newcommand{\rrb}{\right\rbrace}
\newcommand{\C}{\mathbf{C}}

\newtheorem{theorem}{Theorem}

\newtheoremrm{rem}{Remark}

\theoremstyle{definition}

\mathrmdef{Hom}
\mathbfdef{Set}

\begin{document}

\maketitle
\begin{abstract}
A PROB is a ``product and braid" category. Such categories can be used to encode the structure borne by an object in a braided monoidal category. In this paper we provide PROBs whose categories of algebras in a braided monoidal category are equivalent to the categories of monoids and comonoids using the category associated to the braid crossed simplicial group of Fiedorowicz and Loday. We show that PROBs can be composed by generalizing the machinery introduced by Lack for PROPs. We use this to define a PROB for bimonoids in a braided monoidal category as a composite of the PROBs for monoids and comonoids.
\end{abstract}

\section*{Introduction}
\label{sec-Introduction}
A \emph{PROB} is a ``product and braid" category. Such categories are used to encode the structure borne by an object in a braided monoidal category. A PROB is the braided monoidal analogue of a \emph{PROP} in the symmetric monoidal setting and a \emph{Lawvere theory} in the cartesian monoidal setting. 

In this paper we present PROBs whose categories of algebras in a braided monoidal category are equivalent to the categories of monoids, comonoids and bimonoids. For monoids and comonoids the PROBs are closely related to the category associated to the braid crossed simplicial group of Fiedorowicz and Loday \cite[3.7]{FL}. We demonstrate that Lack's methods for composing PROPs \cite{Lack} generalize to give a notion of composing PROBs. We use these methods to form a composite PROB from the PROBs for monoids and comonoids and demonstrate that the category of algebras in a braided monoidal category for this composite is equivalent to the category of bimonoids. The results of this paper can be seen both as an extension of the theory introduced in \cite{Lack} to the setting of braided monoidal categories and an extension of using the structure inherent in a crossed simplicial group to categorify objects in a symmetric monoidal category as studied in \cite{Pir-PROP}, \cite{Lack} and \cite{ifas}. In particular our main theorem, Theorem \ref{bimon-thm}, can be seen as a braided monoidal analogue of \cite[Theorem 5.2]{Pir-PROP} and \cite[5.9]{Lack}.

The paper is organized as follows. In Section \ref{PROB-sec} we recall the definitions of PROs, PROBs and PROPs. In Section \ref{EH-sec} we give a version of the Eckmann-Hilton argument for PROs, PROBs and PROPs. In Section \ref{EG-sec} we give examples of PROs and PROBs. We recall the notion of a distributive law of PROs and use it to construct a PRO, with a canonical PROB structure, from the PRO of finite ordinals and the PROB of braid groups. This PROB is denoted $\mathbb{D}\otimes \mathbb{B}$. In Section \ref{YB-sec} we recall the connection between the PROB of braid groups and Yang-Baxter operators. In Section \ref{comp-sec} we provide analogues to Lack's results on composing PROPs. We construct the notion of a distributive law for PROBs and a composite of PROBs. We provide results describing the structure of an algebra for a composite PROB in terms of algebras for the two factors of the composite. In Section \ref{mon-sec} we prove that the category of algebras for $\mathbb{D}\otimes \mathbb{B}$ in a braided monoidal category is equivalent to the category of monoids and that the category of algebras for its opposite PROB is equivalent to the category of comonoids. In Section \ref{bimon-sec} we define a distributive law between the PROB $\mathbb{D}\otimes \mathbb{B}$ and its opposite to obtain a composite PROB. We prove that the category of algebras in a braided monoidal category for the composite PROB is equivalent to the category of bimonoids.

\subsection*{Acknowledgements}
I am very grateful to Callum Reader, James Cranch and James Brotherston for their helpful conversations and suggestions. I would like to thank Ross Street for alerting me to a revision in the paper \cite{Lack}. I would like to thank the referee for their helpful comments and interesting questions.

\section{PROBs}
\label{PROB-sec}

\definition
\label{sets-defn}
For $n\geqslant 1$ we define $\ul{n}$ to be the set $\llb 1,\dotsc ,n\rrb$. We define $\ul{0}=\emptyset$.
\enddefinition

\definition
A \emph{PRO} is a strict monoidal category whose objects are the sets $\ul{n}$ for $n\geqslant 0$ and whose tensor product is given by addition. 
\enddefinition

\definition
A \emph{PROB} is a braided strict monoidal category whose objects are the sets $\ul{n}$ for $n\geqslant 0$ with tensor product given by addition.
\enddefinition

\definition
A \emph{PROP} is a symmetric strict monoidal category whose objects are the sets $\ul{n}$ for $n\geqslant 0$ with tensor product given by addition.
\enddefinition

Let $\mathbf{MonCat}$, $\mathbf{BrMonCat}$ and $\mathbf{SymMonCat}$ denote the $2$-categories of monoidal categories, braided monoidal categories and symmetric monoidal categories respectively. The objects are monoidal categories, braided monoidal categories and symmetric monoidal categories respectively. The morphisms are monoidal functors, braided monoidal functors and symmetric monoidal functors respectively. The $2$-morphisms in each case are morphisms of monoidal functors, that is, natural transformations that are compatible with the monoidal structure (see \cite[Section 1]{JS} for instance).

\definition
We denote by $\mathbf{PRO}$, $\mathbf{PROB}$ and $\mathbf{PROP}$ the $2$-categories of PROs, PROBs and PROPs respectively. The objects are PROs, PROBs and PROPs respectively. The morphisms are identity-on-objects strict monoidal functors, identity-on-objects braided strict monoidal functors and identity-on-objects symmetric strict monoidal functors respectively. The $2$-morphisms in each case are morphisms of monoidal functors.
\enddefinition

We see that $\mathbf{PRO}$, $\mathbf{PROB}$ and $\mathbf{PROP}$ are sub-$2$-categories of $\mathbf{MonCat}$, $\mathbf{BrMonCat}$ and $\mathbf{SymMonCat}$ respectively. We also observe that $\mathbf{PROP}$ is the full sub-$2$-category of $\mathbf{PROB}$ for which the braidings square to the identity. We can consider the $2$-categories $\mathbf{PRO}$, $\mathbf{PROB}$ and $\mathbf{PROP}$ as categories by forgetting the $2$-morphism structure in the usual way.

\section{The Eckmann-Hilton argument}
\label{EH-sec}
One interpretation of the Eckmann-Hilton argument \cite{EH} states that a monoid in the category of monoids is a commutative monoid. Joyal and Street \cite[Section 5]{JS} generalize this by defining the notion of a multiplication on a monoidal category and proving an equivalence of $2$-categories between $\mathbf{BrMonCat}$ and the category of monoidal categories with multiplication. We will now present a version of their result for PROs and PROBs. 

Let $\mathbf{M}$ be a monoidal category with unit object $I$. Let $1_{\mathbf{M}}$ denote the identity functor on $\mathbf{M}$. Recall from \cite[Section 5]{JS} that a monoidal category $\mathbf{M}$ has a \emph{multiplication} if it is equipped with a monoidal functor $\Phi\colon \mathbf{M}\times \mathbf{M} \rightarrow \mathbf{M}$ together with natural isomorphisms $\rho\colon \Phi\circ \left(1_M,I\right)\rightarrow 1_M$ and $\lambda\colon  \Phi\circ \left(I,1_M\right)\rightarrow 1_M$.

\definition
Let $\mathrm{Mult}\left(\mathbf{PRO}\right)$ denote the $2$-category of PROs with multiplication. Let $\mathrm{Mult}\left(\mathrm{Mult}\left(\mathbf{PRO}\right)\right)$ denote the $2$-category of PROs in $\mathrm{Mult}\left(\mathbf{PRO}\right)$ with multiplication. Let $\mathrm{Mult}\left(\mathbf{PROB}\right)$ denote the $2$-category of PROBs with multiplication.
\enddefinition

\proposition
There are equivalences of $2$-categories
\[\mathbf{PROB}\simeq \mathrm{Mult}\left(\mathbf{PRO}\right)\]
and
\[\mathbf{PROP} \simeq \mathrm{Mult}\left(\mathbf{PROB}\right) \simeq \mathrm{Mult}\left(\mathrm{Mult}\left(\mathbf{PRO}\right)\right).\]
\endproposition
\proof
The result follows from specific instances of Propositions 5.2, 5.3 and 5.4 in \cite{JS}.
\endproof

\section{Algebras, examples and distributive laws}
\label{EG-sec}

\definition
Let $\mathbf{T}$ be a PRO. For a monoidal category $\mathbf{M}$, an \emph{algebra of $\mathbf{T}$ in $\mathbf{M}$} is a strict monoidal functor $\mathbf{T}\rightarrow \mathbf{M}$. 
\enddefinition

\definition
Let $\mathbf{P}$ be a PROB. For a braided monoidal category $\mathbf{B}$, a \emph{$\mathbf{P}$-algebra in $\mathbf{B}$} is a braided strict monoidal functor $\mathbf{P}\rightarrow\mathbf{B}$. We denote the category of $\mathbf{P}$-algebras in $\mathbf{B}$ and natural transformations by $\mathbf{Alg}\left(\mathbf{P},\mathbf{B}\right)$.
\enddefinition

\example
We denote by $\mathbb{D}$ the PRO of finite ordinals and order-preserving maps as in \cite[2.2]{Lack}. For a strict monoidal category $\mathbf{M}$, an algebra of $\mathbb{D}$ in $\mathbf{M}$ is a monoid in $\mathbf{M}$, see \cite[VII 5]{CWM}.
\endexample

\example
Recall the PRO of braid groups, denoted $\mathbb{B}$, from \cite[Example 2.1]{JS}. The set $\mathrm{Hom}_{\mathbb{B}}\left(\ul{n},\ul{m}\right)$ is empty for $n\neq m$ and $\mathrm{Hom}_{\mathbb{B}}\left(\ul{n} , \ul{n}\right)=B_n$, the braid group of $n$ strings. The strict monoidal structure is given by the \emph{addition of braids}. 
\endexample

\remark
Observe that a PROB determines a PRO $\mathbf{T}$ with a morphism of PROs $F\colon \mathbb{B}\rightarrow \mathbf{T}$. The braiding $\ul{m}+\ul{n}\rightarrow \ul{n}+\ul{m}$ in $\mathbf{T}$ is the image under the functor $F$ of the braid which interchanges the first $m$ elements with the last $n$ elements, braiding the former over the latter.

Conversely, such a morphism of PROs arises from a PROB if and only if the resulting maps $\ul{m}+\ul{n}\rightarrow \ul{n}+\ul{m}$ are natural. A morphism of PROBs is a morphism of PROs which commutes with the maps out of $\mathbb{B}$. We therefore see that the category $\mathbf{PROB}$ is a full subcategory of the under-category $\mathbb{B}\downarrow \mathbf{PRO}$.

In particular, we observe that $\mathbb{B}$ has a canonical PROB structure given by the identity functor.
\endremark

\remark
\label{PRO-alg-rem}
Given a PROB $\mathbf{P}$ we will sometimes want to refer to an algebra structure for the underlying PRO. In this case we will refer to a $\mathbf{P}$-PRO-algebra. This will used in Proposition \ref{monoid-prop}.
\endremark

\remark
\label{BMC-rem}
Let $\mathbf{1}$ denote the category with one object and one arrow. By \cite[Proposition 2.2(b)]{JS}, the PROB $\mathbb{B}$ is the free braided strict monoidal category on the category $\mathbf{1}$. It follows from \cite[Corollary 2.4]{JS} that for any braided monoidal category $\mathbf{B}$ there is an equivalence of categories $\mathbf{Alg}\left(\mathbb{B},\mathbf{B}\right) \simeq \mathbf{B}$.
\endremark

We define a further example that arise as a consequence of a distributive law. Let $\mathbb{S}$ and $\mathbb{T}$ be PROs. Recall from \cite[Section 3]{Lack} that a distributive law of PROs $L\colon \mathbb{T}\otimes \mathbb{S}\rightarrow \mathbb{S}\otimes \mathbb{T}$ is defined as follows. Given $\sigma\in \mathrm{Hom}_{\mathbb{S}} \left(\ul{m}, \ul{n}\right)$ and $\tau\in \mathrm{Hom}_{\mathbb{T}}\left(\ul{n},\ul{r}\right)$, one has an object $\sigma_L\tau$ and morphisms $\sigma_{\mathbb{T}}\tau \in \mathrm{Hom}_{\mathbb{T}}\left(\ul{m}, \sigma_L\tau\right)$ and $\sigma_{\mathbb{S}}\tau \in \mathrm{Hom}_{\mathbb{S}}\left( \sigma_L \tau , \ul{r}\right)$, satisfying the equations of \cite[Section 2]{RW} and compatible with the monoidal structure.

\example
\label{DB-eg}
Given a pair $\left(h, \psi\right)$ where $\psi \in \mathrm{Hom}_{\mathbb{D}}\left(\ul{n} , \ul{m}\right)$ and $h\in \mathrm{Hom}_{\mathbb{B}}\left(\ul{m} , \ul{m}\right)$ there is a unique pair $\left(h_{\star}\left(\psi\right),\psi^{\star}(h)\right)$ where $\psi^{\star}(h)\in \mathrm{Hom}_{\mathbb{B}}\left(\ul{n} ,\ul{n}\right)$ and $h_{\star}(\psi)\in \mathrm{Hom}_{\mathbb{D}}\left(\ul{n} , \ul{m}\right)$, as constructed in \cite[3.7]{FL} (see also \cite[Section 4]{DS}). The fact that these assignments satisfy the relations of a distributive law \cite[2.4]{RW} follows from the fact that they satisfy the relations of a crossed simplicial group given in \cite[1.6]{FL} and a routine check shows that they respect the monoidal structures of $\mathbb{D}$ and $\mathbb{B}$. By \cite[3.8]{Lack} we have a PRO $\mathbb{D}\otimes \mathbb{B}$ whose morphisms are the pairs of the form $\left(\phi , g\right)$ where $g\in \mathrm{Hom}_{\mathbb{B}}\left(\ul{n} , \ul{n}\right)$ and $\phi\in \mathrm{Hom}_{\mathbb{D}}\left(\ul{n} ,\ul{m}\right)$ with composition defined via the distributive law. It has a canonical PROB structure induced from $\mathbb{B}$.
\endexample

\section{Yang-Baxter operators}
\label{YB-sec}

The PROB $\mathbb{B}$ is closely related to the study of Yang-Baxter operators. Yang-Baxter operators appear in a variety of settings; notably in the the study of link invariants \cite{turaev} and quantum groups \cite{drinfeld}, \cite{manin}. 

Following \cite[Definition 2.4]{JS}, a \emph{Yang-Baxter operator} on a functor $T\colon \C \rightarrow \mathbf{M}$, where $\mathbf{M}$ is a monoidal category, is a natural family of isomorphisms
\[y_{C_1,C_2}\colon T\left(C_1\right) \otimes T\left(C_2\right) \rightarrow T\left(C_2\right) \otimes T\left(C_1\right)\]
subject to a compatibility condition known as the \emph{Yang-Baxter equation}. For example, any functor whose target category is braided monoidal comes with a Yang-Baxter operator given by the braiding. 

Recall the category $\mathbf{1}$ from Remark \ref{BMC-rem}. Let $\mathbf{M}$ be a strict monoidal category. A Yang-Baxter operator on a functor $T\colon \mathbf{1}\rightarrow \mathbf{M}$ consists of an object $X\in \mathbf{M}$ together with an isomorphism $y\colon X\otimes X\rightarrow X\otimes X$ satisfying the Yang-Baxter equation. Denote by $\mathbf{YB}\left(\mathbf{M}\right)$ the category of Yang-Baxter operators in $\mathbf{M}$. The objects are pairs $\left(X, y\right)$ where $X$ is an object of $\mathbf{M}$ and $y$ is a Yang-Baxter operator on $X$. A morphism $\left(X_1, y_1\right) \rightarrow \left(X_2 , y_2\right)$ in $\mathbf{YB}\left(\mathbf{M}\right)$ is a morphism $X_1\rightarrow X_2$ in $\mathbf{M}$ compatible with the Yang-Baxter structure.

Using the terminology of Remark \ref{PRO-alg-rem}, \cite[Proposition 2.2(a)]{JS} tells us that the category of $\mathbb{B}$-PRO-algebras in $\mathbf{M}$ is equivalent to the category $\mathbf{YB}\left(\mathbf{M}\right)$ of Yang-Baxter operators in $\mathbf{M}$.

\section{Composing PROBs}
\label{comp-sec}
We follow the methods of \cite[Section 4]{Lack} to show that one can compose PROBs analogously to the way in which one can compose PROPs. As in the case for PROPs \cite[4.1]{Lack}, given two PROBs $\mathbb{S}$ and $\mathbb{T}$, together with a distributive law of their underlying PROs, we obtain a composite PRO $\mathbb{S}\otimes \mathbb{T}$. This composite PRO comes with two PROB structures coming from the functors $\mathbb{B}\rightarrow \mathbb{S}\rightarrow \mathbb{S}\otimes \mathbb{T}$ and $\mathbb{B}\rightarrow \mathbb{T}\rightarrow \mathbb{S}\otimes \mathbb{T}$. Following \cite[4.3]{Lack}, replacing the PROP of finite sets and permutations with the PROB of braid groups, we take the coequalizer of these functors in $\mathbf{PRO}$, the category of PROs and strict monoidal functors. This allows us to define the notion of a distributive law for PROBs, that is, the data required to give $\mathbb{S}\otimes \mathbb{T}$ a canonical PROB structure. 

We provide analogues of \cite[4.6, 4.7]{Lack}, Theorem \ref{comp-prob-thm} and Proposition \ref{comp-prop} respectively. The former is a result providing conditions under which a PROB may be expressed as a composite via a distributive law. The latter is a result describing the structure of an algebra for a composite PROB in terms of the algebras of the factors.

\definition
Let $\mathcal{N}=\mathbf{Span}\left(\mathbf{Mon}\right)\left(\mathbb{N},\mathbb{N}\right)$ denote the monoidal category of spans of monoids from $\mathbb{N}$ to $\mathbb{N}$, with tensor product given by the composition of spans.
\enddefinition

\remark
As noted in \cite[4.3]{Lack} the category $\mathcal{N}$ has colimits, preserved by tensoring on either side. Furthermore, as noted in \cite[3.5]{Lack} the category $\mathbf{PRO}$ is the category $\mathbf{Mon}\left(\mathcal{N}\right)$ of monoids in $\mathcal{N}$. In particular, the PRO of braid groups $\mathbb{B}$ is a monoid in $\mathcal{N}$. 

\definition 
A \emph{$\mathbb{B}$-bimodule}, $M$, in $\mathbf{PRO}$ consists of the following data:
\begin{itemize}
\item since $M$ is a PRO we have sets of morphisms $\mathrm{Hom}_M\left(\ul{m},\ul{n}\right)$ for all $m$ and $n$ in $\mathbb{N}$;
\item we have left and right actions of the braid groups on these sets;
\item there is a unital, associative operation \[\amalg\colon \mathrm{Hom}_M\left(\ul{m_1},\ul{n_1}\right) \times \mathrm{Hom}_M\left(\ul{m_2},\ul{n_2}\right) \rightarrow \mathrm{Hom}_M\left(\ul{m_1+n_1},\ul{m_2+n_2}\right)\] induced from the monoidal structure of $M$;
\item this operation is compatible with the actions of the braid groups in the sense that \[\left(\pi_1\amalg \pi_2\right)\left(f_1\amalg f_2\right)\left(\sigma_1\amalg \sigma_2\right)= \pi_1f_1\sigma_1 \amalg \pi_2f_2\sigma_2\] for $f_i\in \mathrm{Hom}_M\left(\ul{m_i},\ul{n_i}\right)$, $\pi_i\in B_{n_i}$ and $\sigma_i\in B_{m_i}$.
\end{itemize}
Let $\mathcal{B}^{\prime}$ denote the category of \emph{$\mathbb{B}$-bimodules}.
\enddefinition

\definition
\label{coeq-defn}
Let $M$ and $N$ be $\mathbb{B}$-bimodules. Let $\lambda$ denote the left action of $\mathbb{B}$ on $N$ and let $\rho$ denote the right action of $\mathbb{B}$ on $M$. The tensor product of $\mathbb{B}$-bimodules, $M\otimes_{\mathbb{B}} N$, is defined by the coequalizer
\begin{center}
\begin{tikzcd}
M\otimes \mathbb{B} \otimes N \ar[r,shift left=.75ex,"{\rho \otimes N}"]
  \ar[r,shift right=.75ex,swap,"{M\otimes \lambda}"]
&
M \otimes N \ar[r] 
&
M\otimes_{\mathbb{B}} N
\end{tikzcd}
\end{center}
in $\mathbf{PRO}$.
\enddefinition

\remark
As in the case for PROPs, the tensor product $-\otimes_{\mathbb{B}}-$ endows the category $\mathcal{B}^{\prime}$ with the structure of a monoidal category.
\endremark

Analogously to \cite[4.3]{Lack}, the category $\mathbf{Mon}\left(\mathcal{B}^{\prime}\right)$ of monoids in $\mathcal{B}^{\prime}$ is equivalent to the under-category $\mathbb{B}\downarrow \mathbf{Mon}\left(\mathcal{N}\right)$, that is, it is equivalent to the under-category $\mathbb{B}\downarrow \mathbf{PRO}$.

\definition
Consider the braids $\ul{m_1}+\ul{m_2}\rightarrow \ul{m_2}+\ul{m_1}$ and $\ul{n_1}+\ul{n_2}\rightarrow \ul{n_2}+\ul{n_1}$ which interchange the first $m_1$ (respectively $n_1$) elements with the last $m_2$ (respectively $n_2$) elements, braiding the former over the latter.

Let $\mathcal{B}$ denote the full subcategory of $\mathcal{B}^{\prime}$ which consists of those $\mathbb{B}$-bimodules $M$ for which acting on $f_1\amalg f_2\in \mathrm{Hom}_M\left(\ul{m_1+m_2}, \ul{n_1+n_2}\right)$ by these braids gives $f_2\amalg f_1$.
\enddefinition

We observe that the monoids in the category $\mathcal{B}$ are precisely the PROBs.

\definition
We define a distributive law of PROBs to be a distributive law of monoids in $\mathcal{B}$: that is, a morphism $\mathbb{T}\otimes_{\mathbb{B}}\mathbb{S}\rightarrow \mathbb{S}\otimes_{\mathbb{B}} \mathbb{T}$ of $\mathbb{B}$-bimodules satisfying the equations of \cite[Section 2]{RW}.
\enddefinition

\remark
Similarly to the case for PROPs \cite[4.4]{Lack} we can express this construction as a monad on an object in a $2$-category in the sense of \cite[Section 1]{FTM}. Let $\mathbf{Mon}$ be the category of monoids in sets and monoid homomorphisms. Let $\mathbf{Prof}\left(\mathbf{Mon}\right)$ denote the $2$-category whose objects are categories internal to $\mathbf{Mon}$, whose morphisms are internal profunctors and whose $2$-morphisms are natural transformations. The monoidal category $\mathbb{B}$ is an object of $\mathbf{Prof}\left(\mathbf{Mon}\right)$ and a $\mathbb{B}$-bimodule is a morphism in $\mathbf{Prof}\left(\mathbf{Mon}\right)$ from $\mathbb{B}$ to $\mathbb{B}$. That is, the category $\mathcal{B}^{\prime}$ can be identified with $\mathbf{Prof}\left(\mathbf{Mon}\right)\left(\mathbb{B},\mathbb{B}\right)$. Therefore a PROB is precisely a monad in $\mathbf{Prof}\left(\mathbf{Mon}\right)$ on the object $\mathbb{B}$ whose underlying bimodule lies not only in $\mathcal{B}^{\prime}$ but also in $\mathcal{B}$.

It is also worth remarking that similar constructions apply in the theory of distributive laws for Lawvere theories by work of Cheng \cite[Section 5]{Cheng}.
\endremark

\theorem
\label{comp-prob-thm}
Let $\mathbb{R}$ be a PROB. Let $\mathbb{S}$ and $\mathbb{T}$ be subcategories of $\mathbb{R}$ containing all the objects, all the braidings and closed under tensoring. Suppose that every morphism $\rho$ in $\mathbb{R}$ can be written as a composite $\rho=\sigma \circ \tau$ where $\tau$ is a morphism in $\mathbb{T}$ and $\sigma$ is a morphism in $\mathbb{S}$. Furthermore suppose that if $\rho= \sigma^{\prime} \circ \tau^{\prime}$ is another such representation then there is a braid $\pi$ such that $\tau = \tau^{\prime}\circ \pi$ and $\pi\circ \sigma = \sigma^{\prime}$. Then $\mathbb{R}$ is the composite of $\mathbb{S}$ and $\mathbb{T}$ via a distributive law $L\colon \mathbb{T}\otimes_{\mathbb{B}}\mathbb{S}\rightarrow \mathbb{S}\otimes_{\mathbb{B}} \mathbb{T}$.
\endtheorem
\proof
This is analogous to \cite[4.6]{Lack}.
\endproof

\proposition
\label{comp-prop}
If $\mathbb{S}\otimes_{\mathbb{B}} \mathbb{T}$ is a composite PROB induced by a distributive law then an $\mathbb{S}\otimes_{\mathbb{B}}\mathbb{T}$-algebra structure on an object $A$ in a braided monoidal category $\mathbf{B}$ consists of an $\mathbb{S}$-algebra structure and a $\mathbb{T}$-algebra structure subject to the condition that 
\begin{center}
\begin{tikzcd}
A^{\otimes m}\arrow[r, "\sigma "]\arrow[d, "{\sigma_{\mathbb{T}}\tau}",swap] & A^{\otimes n}\arrow[d, "\tau "]\\
A^{\otimes q}\arrow[r, "{\sigma_{\mathbb{S}}\tau}",swap] & A^{\otimes r}
\end{tikzcd}
\end{center}
where $\sigma\in \mathrm{Hom}_{\mathbb{S}}\left(\ul{m} , \ul{n}\right)$, $\tau\in \mathrm{Hom}_{\mathbb{T}}\left(\ul{n} , \ul{r}\right)$ and $q=\sigma_L\tau$.

Furthermore, if objects $A$ and $B$ in $\mathbf{B}$ have $\mathbb{S}\otimes_{\mathbb{B}} \mathbb{T}$-algebra structures then a morphism $f\in \mathrm{Hom}_{\mathbf{B}}\left(A,B\right)$ is a morphism of $\mathbb{S}\otimes_{\mathbb{B}}\mathbb{T}$-algebras if and only if it is a morphism of $\mathbb{S}$-algebras and a morphism of $\mathbb{T}$-algebras.
\endproposition
\proof
This is analogous to \cite[3.9--3.12]{Lack}. See also \cite[4.7, 4.8]{Lack}.
\endproof

As noted in Section \ref{PROB-sec}, $\mathbf{PROP}$ is the full sub-$2$-category of $\mathbf{PROB}$ for which the braidings square to the identity. Therefore, every PROP can be thought of as a PROB with the additional condition that the braiding squares to the identity. The following proposition shows that our notion of composing PROBs is compatible with the composition of PROPs.

\proposition
Let $P_1$ and $P_2$ be PROPs. There is an equality of PROPs $P_1\otimes_{\mathbb{P}} P_2 = P_1\otimes_{\mathbb{B}} P_2$.
\endproposition
\proof
$P_1$ and $P_2$ are PROBs such that the braiding squares to the identity. In other words the braidings factor through the symmetric groups via the functor $\mathbb{B}\rightarrow \mathbb{P}$ which sends a braid to its underlying permutation. This means that the identifications in the coequalizer
\begin{center}
\begin{tikzcd}
P_1\otimes \mathbb{B} \otimes P_2 \ar[r,shift left=.75ex,"{\rho \otimes P_2}"]
  \ar[r,shift right=.75ex,swap,"{P_1\otimes \lambda}"]
&
P_1 \otimes P_2 \ar[r] 
&
P_1\otimes_{\mathbb{B}} P_2
\end{tikzcd}
\end{center}
of Definition \ref{coeq-defn} coincide with the identifications in the coequalizer
\begin{center}
\begin{tikzcd}
P_1\otimes \mathbb{P} \otimes P_2 \ar[r,shift left=.75ex,"{\rho \otimes P_2}"]
  \ar[r,shift right=.75ex,swap,"{P_1\otimes \lambda}"]
&
P_1 \otimes P_2 \ar[r] 
&
P_1\otimes_{\mathbb{P}} P_2
\end{tikzcd}
\end{center}
of \cite[4.3]{Lack}, from which the result follows.
\endproof

\section{PROBs for monoids and comonoids}
\label{mon-sec}
We show that the category of algebras in a braided monoidal category $\mathbf{B}$ for the PROB $\mathbb{D}\otimes \mathbb{B}$ of Example \ref{DB-eg} is equivalent to the category of monoids in $\mathbf{B}$. As a corollary we show that the category of algebras in $\mathbf{B}$ for the PROB $\left(\mathbb{D}\otimes \mathbb{B}\right)^{op}$ is equivalent to the category of comonoids in $\mathbf{B}$.

\definition
We denote the category of monoids in a braided monoidal category $\mathbf{B}$ by $\mathbf{Mon}\left(\mathbf{B}\right)$. We denote the category of comonoids in $\mathbf{B}$ by $\mathbf{Comon}\left(\mathbf{B}\right) = \mathbf{Mon}\left(\mathbf{B}^{op}\right)^{op}$. Let $\mathbf{Bimon}\left(\mathbf{B}\right) = \mathbf{Mon}\left(\mathbf{Comon}\left(\mathbf{B}\right)\right) = \mathbf{Comon}\left(\mathbf{Mon}\left(\mathbf{B}\right)\right)$ denote the category of bimonoids in $\mathbf{B}$.
\enddefinition

\proposition
\label{monoid-prop}
Let $\mathbf{B}$ be a braided monoidal category. There is an equivalence of categories $\mathbf{Alg}\left( \mathbb{D}\otimes \mathbb{B} , \mathbf{B}\right) \simeq \mathbf{Mon}\left(\mathbf{B}\right)$.
\endproposition
\proof
Recall the terminology of Remark \ref{PRO-alg-rem}. By \cite[3.10]{Lack}, a $\left(\mathbb{D}\otimes \mathbb{B}\right)$-PRO-algebra structure on an object $M$ of $\mathbf{B}$ consists of a $\mathbb{D}$-algebra structure and a $\mathbb{B}$-algebra structure subject to a compatibility condition. A $\mathbb{D}$-algebra structure is a monoid structure. A $\mathbb{B}$-algebra is an object $M$ together with an isomorphism $M^{\otimes n}\rightarrow M^{\otimes n}$ for each element of the braid group $B_n$. Arguing analogously to \cite[5.5]{Lack} a $\left(\mathbb{D}\otimes \mathbb{B}\right)$-PRO-algebra structure is a $\left(\mathbb{D}\otimes \mathbb{B}\right)$-algebra structure if and only if the only isomorphisms $M^{\otimes n}\rightarrow M^{\otimes n}$ are those induced from the braidings. The compatibility condition follows from the naturality of the braidings. Finally, a morphism in $\mathbf{B}$ is a map of monoids if and only if it respects the $\mathbb{D}$-algebra structure and the $\mathbb{B}$-algebra structure. By Proposition \ref{comp-prop}, this is true if and only if it respects the $\left(\mathbb{D}\otimes \mathbb{B}\right)$-algebra structure.
\endproof

\corollary
\label{comon-cor}
Let $\mathbf{B}$ be a braided monoidal category. There is an equivalence of categories $\mathbf{Alg}\left( \left(\mathbb{D}\otimes \mathbb{B}\right)^{op} , \mathbf{B}\right) \simeq \mathbf{Comon}\left(\mathbf{B}\right)$.
\endcorollary
\proof
Using Proposition \ref{monoid-prop} we observe that
\[\mathbf{Comon}\left(\mathbf{B}\right)=\mathbf{Mon}\left(\mathbf{B}^{op}\right)^{op}\simeq \mathbf{Alg}\left( \left(\mathbb{D}\otimes \mathbb{B}\right) , \mathbf{B}^{op}\right)^{op}=\mathbf{Alg}\left( \left(\mathbb{D}\otimes \mathbb{B}\right)^{op} , \mathbf{B}\right)\]
as required.
\endproof

\section{The PROB for bimonoids}
\label{bimon-sec}
In this section we will define a distributive law of PROBs between $\mathbb{D}\otimes \mathbb{B}$, the PROB for monoids, and $\left(\mathbb{D}\otimes \mathbb{B}\right)^{op}$, the PROB for comonoids and prove that the algebras for the composite are bimonoids.

The distributive law required takes the form of a map of $\mathbb{B}$-bimodules 
\[\left(\mathbb{D}\otimes\mathbb{B}\right)^{op}\otimes_{\mathbb{B}}\left(\mathbb{D}\otimes \mathbb{B}\right)\rightarrow \left(\mathbb{D}\otimes \mathbb{B}\right)\otimes_{\mathbb{B}}\left(\mathbb{D}\otimes \mathbb{B}\right)^{op}\]
from the PROB of equivalences classes of cospans in $\mathbb{D}\otimes\mathbb{B}$ to the category of equivalence classes of spans in $\mathbb{D}\otimes\mathbb{B}$, subject to the conditions of \cite[Section 2]{RW}.

The procedure for doing this is analogous to the methods employed by Pirashvili in the construction of the double category $\mathcal{F}(as)$ \cite[Section 4]{Pir-PROP} and Lack \cite[5.9]{Lack} in the symmetric monoidal case. Given an equivalence class of cospans 
\begin{center}
\begin{tikzcd}
\ul{p}\arrow[rr, "{\left(\phi ,g\right)}"]&&\ul{n}&&\ul{q}\arrow[ll, "{\left(\psi ,h\right)}",swap] 
\end{tikzcd}
\end{center}
in $\mathbb{D}\otimes \mathbb{B}$, we will take the pullback in the category of finite sets and use the data of the distributive law for $\mathbb{D}\otimes \mathbb{B}$ to give a unique lift of the pullback to the category of equivalence classes of spans in $\mathbb{D}\otimes \mathbb{B}$.

Recall that a morphism in $\mathbb{D}\otimes \mathbb{B}$ is a unique pair $\left(\phi , g\right)$ where $g$ is an element of a braid group and $\phi$ is an order-preserving map, with composition defined by the distributive law of \cite[3.7]{FL}. Recall that the elements of the braid groups are the isomorphisms in $\mathbb{D}\otimes \mathbb{B}$. Furthermore, using composition and the disjoint union, all morphisms in $\mathbb{D}$ have an expression in terms of the unique morphisms $m\in \mathrm{Hom}_{\mathbb{D}}\left(\ul{2},\ul{1}\right)$ and $u\in \mathrm{Hom}_{\mathbb{D}}\left(\ul{0},\ul{1}\right)$ with finitely many terms. 

In the following definition, Points \ref{braid1} and \ref{mixed} give the necessary data when at least one of the morphisms in the cospan is an isomorphism. Point \ref{ord-pres1} gives the necessary data when both morphisms in the cospan are order-preserving maps. By combining the properties of pullback diagrams, namely compatibility with the disjoint union and composition of set maps, and the fact that every morphism in $\mathbb{D}$ has a finite expression in terms of $m$ and $u$, it suffices to give the assignment on the given cospans.

We observe that the data of the following definition satisfies the conditions of a distributive law by construction. It is straightforward to see that identities are preserved in the sense of \cite{RW}. Furthermore, the assignments are compatible with composition in $\mathbb{D}\otimes \mathbb{B}$ and the monoidal structure since any morphism in $\mathbb{D}\otimes \mathbb{B}$ can be written in terms of the braid groups and the morphisms $m$ and $u$ using coproducts and the distributive law $\mathbb{D}\otimes \mathbb{B}$ found in \cite[3.7]{FL}.

\definition
\label{dist-law-defn}
We define a distributive law of PROBs
\[\left(\mathbb{D}\otimes\mathbb{B}\right)^{op}\otimes_{\mathbb{B}}\left(\mathbb{D}\otimes \mathbb{B}\right)\rightarrow \left(\mathbb{D}\otimes \mathbb{B}\right)\otimes_{\mathbb{B}}\left(\mathbb{D}\otimes \mathbb{B}\right)^{op}\]
to be determined as follows
\begin{enumerate}
\item \label{braid1} For $g$, $h\in B_n$
\begin{small}
\begin{center}
\begin{tikzcd}
\ul{n}\arrow[rr, "g"]&&\ul{n}&&\ul{n}\arrow[ll, "h",swap] & \mapsto & \ul{n}&&\ul{n}\arrow[ll,"{g^{-1}}",swap]\arrow[rr, "{h^{-1}}"]&&\ul{n}
\end{tikzcd}
\end{center}
\end{small}
\item \label{mixed} For $\phi \in \mathrm{Hom}_{\mathbb{D}}\left(\ul{n},\ul{m}\right)$ and $g\in B_m$
\begin{small}
\begin{center}
\begin{tikzcd}
\ul{n}\arrow[rr, "{\phi}"]&&\ul{m}&&\ul{m}\arrow[ll, "g",swap] & \mapsto  & \ul{n}&&\ul{n}\arrow[ll,"{g_{\star}(\phi)}",swap]\arrow[rr, "{\phi^{\star}(g)}"]&&\ul{m}\\
\ul{m}\arrow[rr, "{g}"]&&\ul{m}&&\ul{n}\arrow[ll, "{\phi}",swap] & \mapsto  & \ul{m}&&\ul{n}\arrow[ll,"{\phi^{\star}(g)}",swap]\arrow[rr, "{g_{\star}(\phi)}"]&&\ul{n}
\end{tikzcd}
\end{center}
\end{small}
where $g_{\star}(\phi)$ and $\phi^{\star}(g)$ are the maps determined by the distributive law for $\mathbb{D}\otimes \mathbb{B}$.
\item \label{ord-pres1} For $m\in \mathrm{Hom}_{\mathbb{D}}\left(\ul{2},\ul{1}\right)$ and $u\in \mathrm{Hom}_{\mathbb{D}}\left(\ul{0},\ul{1}\right)$
\begin{small}
\begin{center}
\begin{tikzcd}
\ul{2}\arrow[rr, "m"]&&\ul{1}&&\ul{2}\arrow[ll, "m",swap] & \mapsto  & \ul{2}&&\ul{4}\arrow[ll,"{\left(m+m\right)\circ \sigma_{2,3}}",swap]\arrow[rr, "{m+m}"]&&\ul{2}\\
\ul{0}\arrow[rr, "u"]&&\ul{1}&&\ul{0}\arrow[ll, "u",swap] & \mapsto & \ul{0}&&\ul{0}\arrow[ll,"{id_0}",swap]\arrow[rr, "{id_0}"]&&\ul{0}\\
\ul{2}\arrow[rr, "m"]&&\ul{1}&&\ul{0}\arrow[ll, "u",swap] & \mapsto & \ul{2}&&\ul{0}\arrow[ll,"{u+u}",swap]\arrow[rr, "{id_0}"]&&\ul{0}\\
\ul{0}\arrow[rr, "u"]&&\ul{1}&&\ul{2}\arrow[ll, "m",swap] & \mapsto & \ul{0}&&\ul{0}\arrow[ll,"{id_0}",swap]\arrow[rr, "{u+u}"]&&\ul{2}
\end{tikzcd}
\end{center}
\end{small}
where $\sigma_{2,3}$ is notation for the braid $id_1\amalg \sigma \amalg id_1 \in B_4$ and $\sigma\in B_2$ is the braid which swaps the two elements, braiding the first over the second.
\end{enumerate}
\enddefinition

We will now prove the main theorem, which tells us that the algebras for the composite PROB are bimonoids. Firstly we can tidy up our notation. Since a group is isomorphic to its opposite we have isomorphisms of PROBs
\[\left(\mathbb{D}\otimes \mathbb{B}\right)\otimes_{\mathbb{B}}\left(\mathbb{D}\otimes \mathbb{B}\right)^{op}\cong \left(\mathbb{D}\otimes \mathbb{B}\right)\otimes_{\mathbb{B}}\left(\mathbb{B}\otimes \mathbb{D}^{op}\right)\cong \mathbb{D}\otimes\mathbb{B}\otimes \mathbb{D}^{op}.\]

\theorem
\label{bimon-thm}
Let $\mathbf{B}$ be a braided monoidal category and let $Q=\mathbb{D}\otimes \mathbb{B}\otimes\mathbb{D}^{op}$. There is an equivalence of categories $\mathbf{Alg}\left(Q , \mathbf{B}\right) \simeq \mathbf{Bimon}\left(\mathbf{B}\right)$.
\endtheorem
\proof
By Proposition \ref{comp-prop}, an algebra for $Q$ in $\mathbf{B}$ consists of an object $M$ with a $\left(\mathbb{D}\otimes \mathbb{B}\right)$-algebra structure and a $\left(\mathbb{D}\otimes \mathbb{B}\right)^{op}$-algebra structure subject to the compatibility condition arising from the distributive law. A $\left(\mathbb{D}\otimes \mathbb{B}\right)$-algebra structure is a monoid structure and a $\left(\mathbb{D}\otimes \mathbb{B}\right)^{op}$-algebra structure is a comonoid structure. The compatibility conditions arising from the distributive law of Definition \ref{dist-law-defn} are precisely those requiring $M$ to be a bimonoid.

Finally we observe that a morphism in $\mathbf{B}$ is a morphism of bimonoids if and only if it preserves the $\left(\mathbb{D}\otimes \mathbb{B}\right)$-algebra structure and the $\left(\mathbb{D}\otimes \mathbb{B}\right)^{op}$-algebra structure. By Proposition \ref{comp-prop} this is true if and only if it preserves the $Q$-algebra structure.
\endproof

\end{document}